# An upper bound for the clique number using clique ceiling numbers

R. Dharmarajan[1] and D. Ramachandran[2]

[1]Niels Abel Foundation, Palakkad 678011, Kerala, India.

[1]mathnafrd@gmail.com and [2]claudebergedr@gmail.com

[1]Corresponding author

**Abstract**

In this article we present the idea of clique ceiling numbers of the vertices of a given graph that has a universal vertex. We follow up with a polynomial-time algorithm to compute an upper bound for the clique number of such a graph using clique ceiling numbers. We compare this algorithm with some upper bound formulas for the clique number.

## 1. Introduction

The maximum clique problem (MCP) asks for a clique of the largest possible size in a given graph $G$. Such a clique is called a *maximum clique* of $G$, and the size of any maximum clique of $G$ is the clique number (denoted by $\omega(G)$) of $G$. The MCP is an NP-complete problem [9]. Since the problem is also NP-hard [6], no polynomial-time exact algorithm to solve it is expected to be developed. Nevertheless, it is worthwhile to attempt algorithms for the MCP because the problem has important applications in domains such as social networking, bioinformatics, document clustering, computer vision, image processing and pattern recognition [1,5,11,12].

For an account of upper and lower bounds on $\omega$, we refer the reader to [2,3,4,13]. In this paper we first introduce the idea of clique ceiling numbers. Then we outline an algorithm (named *ACCN*) to compute an upper bound for $\omega$ using clique ceiling numbers. Then we prove the polynomial-time efficiency of the *ACCN* by analysis.

The rest of this paper is organized in sections 2 through 9. Section 2 gives graph theoretical definitions and notation that are relevant to this paper. Section 3 is devoted to the theoretical ideas that will form a base to build the proposed algorithm (*ACCN*) on. Section 4 outlines the *ACCN* in pseudo-code style. In this section we also prove (i) the *ACCN* is feasible and (ii) it is independent of the ordering of the vertices of the input graph. In section 5 we prove that the positive integer returned by the *ACCN* is an upper bound for the clique number of the input graph and then show that the *ACCN* returns precisely the clique number for split graphs. Section 6 analyses the time complexity of the *ACCN*. In section 7 we give a comparison of performances of the *ACCN* with those of three formulas given in [2]. Here we also list merits and limitations of the *ACCN*. An example of computing the clique ceiling numbers is given in section 8. Section 9 is the concluding section.

## 2. Definitions and notation

For the basic definitions and notation used in this article, please see [8]. Throughout this paper, the term graph will mean a simple undirected loop-free graph. If $G = (V, E)$ is a graph, then the expressions $x \in V$ and $x \in G$ will both mean $x$ is a vertex of $G$. Similarly, both $\{x, y\} \in E$ and $\{x, y\} \in G$ will mean $\{x, y\}$ is an edge of $G$. For the rest of this section, $G = (V, E)$ is assumed.

A vertex $u \in G$ is a *universal vertex* of $G$ if $u$ is adjacent to every element of $V - \{u\}$. For $x \in V$, the set $N(x)$ consisting of all the neighbours of $x$ in $G$ is the *neighbourhood* of $x$ in $G$. The *degree* of $x$ in $G$ is denoted by $dx$ or by $dx(G)$, and is the number of vertices of $G$ that are neighbours of $x$ – i.e. $dx = |N(x)|$. A vertex $y$ of $G$ is *isolated* in $G$ if $dy = 0$. $G$ is *null* if $dx = 0$ for every $x \in G$. The smallest degree occuring in $G$ is denoted by $\delta(G)$.

The *closed neighbourhood* of $x$ is denoted by $C(x)$ and is defined as $C(x) = N(x) \cup \{x\}$. Obviously, $u$ is a universal vertex of $G$ if and only if $C(u) = V$.

Two graphs are vertex-disjoint if their vertex sets are disjoint. Let $G_1 = (V_1, E_1)$ and $G_2 = (V_2, E_2)$ be vertex-disjoint graphs. Let $F$ denote the set of all the edges formed by joining each vertex of $G_1$ to each vertex of $G_2$. The *join* of $G_1$ and $G_2$ is denoted by $G_1 \vee G_2$ and is defined to be the graph $(W, L)$ where $W = V_1 \cup V_2$ and $L = E_1 \cup E_2 \cup F$. In particular, if $V_2 = \{z\}$ (so that $E_2 = \phi$) then $G_1 \vee G_2$ is also denoted by $G_1 \vee z$. Obviously $z$ is a universal vertex of $G_1 \vee z$.

A *subgraph* of $G$ is a graph $J = (W, F)$ such that: (i) $W \subset V$, (ii) $F \subset E$ and (iii) each edge in $J$ has the same end points in $J$ as in $G$. $J$ is a *proper subgraph* of $G$ if either $W \neq V$ or $E \neq F$. If $A \subset V$ then the *subgraph induced* by $A$ is the subgraph $G[A] = (A, E[A])$ where $E[A]$ is the set of all those edges $\{x, y\} \in E$ such that $x \in A$ and $y \in A$. In particular, if $a \in V$ then the subgraph induced by $V - \{a\}$ will be denoted by $G - a$.

$G$ is *bipartite* if there is a partition $V = A \cup B$ such that each edge of $G$ has one end in $A$ and the other end in $B$. In this case, $A$ and $B$ are the *partite sets* of $G$ (or, the partite subsets of $V$).

$G$ is *complete* if all of its vertices are pairwise adjacent – i.e., $\{x, y\} \in E$ whenever $x \in V$, $y \in V$ and $x \neq y$. A *clique* of $G$ is a set $M \subset V$ such that $G[M]$ is complete. $M$ is a *maximal clique* of $G$ if (i) $M$ is a clique of $G$ and (ii) $M$ is not a proper subset of any clique of $G$. $M$ is a *maximum clique* of $G$ if (i) $M$ is a clique of $G$ and (ii) $|M| \geq |S|$ for every clique $S$ of $G$.

A graph has a maximum clique though such a clique is not necessarily unique. Obviously, if $M_1$ and $M_2$ are maximum cliques of $G$ then $|M_1| = |M_2|$. If $M$ is a maximum clique of $G$ then the positive integer $|M|$ is the *clique number* ($\omega(G)$) of $G$. If $G$ is null then $\omega(G) = 1$.

The remaining definitions are placed in section 3 since they presuppose propositions therein.

## 3. Preliminaries

Throughout this section, $G = (V, E)$ and $\delta(G) \geq 1$.

**Proposition 3.1.** *Let $x \in G$. For $y \in C(x)$, let $A(y,x) = C(y) \cap C(x)$. Then to each $y \in C(x)$ there corresponds a positive integer $k$ such that there are at least $k - 1$ vertices $z \in A(y, x)$ with $dz(A(y, x)) \geq k - 1$ for each $z$.*

**Proof.** Since $y$ has at least one neighbour in $A(y, x)$, it is immediate that $k = 2$ satisfies the statement.

**Corollary 3.2.** *Let $x \in G$. To each $y \in C(x)$ there corresponds a largest positive integer $k$ that satisfies the statement of proposition 3.1.*

**Proof.** Any positive integer $k$ that satisfies the statement of proposition 3.1 is bounded above by $|G|$ $(= |V|)$.

**Definition.** Let $x \in G$ and $y \in C(x)$. The largest positive integer $k$ that corresponds to $y$ as in corollary 3.2 will be called the *clique ceiling number* (or, *clique ceiling*) of $y$ under $x$, and will be denoted by $c(y{:}x)$.

**Corollary 3.3.** *If $y \in C(x)$ then $c(y{:}x) \leq c(x{:}x)$.*
**Proof.** Let $c(y{:}x) = k$ and $A(y,x) = C(y) \cap C(x)$. Then there exist $z_j \in A(y,x)$ ($j = 1, \ldots, k-1$) such that each $z_j$ has at least $k-1$ neighbours in $A(y,x)$. Consequently $c(x{:}x) \geq k$. ∎

**Proposition 3.4.** *Let $x \in G$. To each $y \in C(x)$ there corresponds a largest positive integer $r$ with the following property: there are at least $r$ vertices $z_i \in C(x) \cap C(y)$ ($i = 1, \ldots, r$), including $x$ and $y$, such that $r \leq c(z_i{:}x)$.*
**Proof.** Certainly $r = 2$ satisfies the statement. Also, any such $r$ is bounded above by $|V|$. ∎

**Definition.** Let $x \in G$. For $y \in C(x)$, the largest positive integer $r$ that corresponds to $y$ as in proposition 3.4 will be called the *revised clique ceiling number* (or, *revised clique ceiling*) of $y$ under $x$, and will be denoted by $c^*(y{:}x)$. The largest $c^*(x{:}x)$ as $x$ runs over $G$ will be called the *clique ceiling number* of $G$, and will be denoted by $c^*(G)$. In symbols, $c^*(G) = \max_{x \in G} \{ c^*(x{:}x) \}$.

**Proposition 3.5.** *Let $u$ be a universal vertex of $G$ and let $x \in G$ be given. Then:*
*(i) $c(y{:}x) \leq c(y{:}u)$ for every $y \in C(x)$ and*
*(ii) $c(x{:}x) \leq c(x{:}u) \leq c(u{:}u)$.*
**Proof.** Obvious, owing to $C(u) = V$.

**Proposition 3.6.** *Let $u$ be a universal vertex of $G$. Then $c^*(u{:}u) \geq c^*(x{:}u) \geq c^*(x{:}x)$ for every $x \in G$.*
**Proof.** Let $x \in G$ be given. Let $c^*(x{:}u) = p$ and $c^*(x{:}x) = q$. By the definition of $c^*$, there exist $p$ vertices $z_i \in C(x)$ ($i = 1, \ldots, p$), including $x$ and $u$, such that $p \leq c(z_i{:}u)$. By dint of these $p$ vertices (and $C(x) \subset C(u)$), we have $c^*(u{:}u) \geq p$.

Next, there exist $q$ vertices $y_i \in C(x)$ ($i = 1, \ldots, q$), including $x$, such that $q \leq c(y_i{:}x)$. By dint of these $q$ vertices (and $C(x) \subset C(u)$), we have $c^*(x{:}u) \geq q$.

**Proposition 3.7.** *Let $G$ be a graph and assume $u \notin G$. If $\omega(G \vee u) \leq r$ then $\omega(G) \leq r - 1$. In particular, if $\omega(G \vee u) = r$ then $\omega(G) = r - 1$.*
**Proof.** Let $\omega(G) = p$ and $M$ be a maximum clique of $G$. Clearly $u \notin M$ and $M \cup \{u\}$ is a clique of $G \vee u$. If $G \vee u$ had a clique of size larger than $p + 1$ then $G \vee u$ would have a clique $A$ such that $u \in A$ and $|A| > p + 1$. This would mean $A - \{u\}$ is a clique of $G$, patently contradicting $\omega(G) = p$. So $\omega(G \vee u) = p + 1$, from which the conclusions follow. ∎

# 4. ACCN – Algorithm for Clique Ceiling Numbers (pseudo-code)

The given instance is $G = (V, E)$, and will be called the *primary instance*. The secondary instance is $G \vee u = (W, F)$ where $u \notin G$, $W = V \cup \{u\}$ and $F$ is the edge set of $G \vee u$. For each $y \in W$, the *ACCN* computes the clique ceiling number of $y$ under $u$, denoted by $c(y:u)$. In the next phase, for each $y \in W - \{u\}$, the *ACCN* computes the revised clique ceiling number of $y$ under $u$, denoted by $c^*(y:u)$. Finally, the *ACCN* computes $c^*(u:u)$ and returns this value.

**OUTLINE OF THE ACCN (Pseudo-code)**
**Input:** Vertex set $V$ and the adjacency matrix of the primary instance $G$.

### Phase 1: Pre-processing

Construct the vertex set and the adjacency matrix of the secondary instance $G \vee u$ and go to phase 2

### Phase 2: Clique ceiling numbers

**for** each $y \in G \vee u$
    $A = C(y)$
    **for** each $z \in A$
      compute $dz(A)$
    **end for**
    compute the largest $k$ such that there are at least $k$ vertices $z$ in $A$ with $dz(A) \geq k$;
    $c(y:u) = k + 1$
**end for** and go to phase 3

### Phase 3: Revised clique ceiling numbers

**for** each $y \in G \vee u$
    $c^*(y:u) = c(y:u)$
**end for**
**for** each $y \in G \vee u$ with $y \neq u$
    compute the largest $p$ such that there are at least $p$ vertices $z \in C(y)$, including $y$ and $u$, such that $p \leq c^*(z:u)$;
    $c^*(y:u) = p$
**end for**
compute the largest $p$ such that there are at least $p$ vertices $z \in G \vee u$, including $u$, such that $p \leq c^*(z:u)$;
$c^*(u:u) = p$ and $c^*(G \vee u) = c^*(u:u)$
Return $c^*(G \vee u)$

-----------------------------------------------------------------------------------------------------------------

**Proposition 4.1.** *The ACCN is feasible – i.e., it terminates in finitely many computations on the instance $G \vee u$.*
**Proof.** Phase 1 requires: (i) adding a single new vertex $u$ to the vertex set of the primary instance $G$ and (ii) expanding the adjacency matrix of $G$ to that of $G \vee u$ by adding one row and one column (equivalently, $2|V|$ new entries). Clearly this phase is feasible.
    Once the control is passed to phase 2, it is never returned to phase 1. Let $G \vee u = (W, F)$ and $|G \vee u| = n$. For each $y \in G \vee u$, $|A| \leq n$. So for each $z \in A$, $dz(A)$ is done in finitely

many computations. So are the computation of $k$ and the assignations of values to the variables $c(y:u)$. Since $|G \vee u| = n$, each 'for' loop in phase 1 terminates in finitely many computations. Hence phase 2 terminates in finitely many computations.

Once the control is passed to phase 3, it is never returned to phase 1 or phase 2. The first 'for' loop in phase 3 terminates in n computations, each of which is an assignation of a positive integer value to a variable. The second 'for' loop in phase 3 terminates in finitely many computations because $|G - u| = n - 1$. Each of these computations involves a logical comparision and an assignation of a positive integer to a variable. The last part of phase 3 terminates in n computations since it runs over each vertex of $C(u)$ (= $W$) and $|C(u)| = n$. Hence phase 3 terminates in finitely many computations. ∎

**Proposition 4.2.** *If the ACCN returns $c^*(G \vee u) = r$ in some ordering of the vertices of $G \vee u$, then it returns $c^*(G \vee u) = r$ in every ordering of the vertices of $G \vee u$. Hence the ACCN is independent of the ordering of the vertices of the secondary input.*

**Proof.** Let two orderings of the vertices of $G \vee u$ be given. Let $x \in G \vee u$ be arbitrary. In these two orderings, $C(x)$ is the same. Then so is $J$, the subgraph (of $G \vee u$) induced by $C(x)$. Then so are $dy(J)$ for every $y \in C(x)$ and (hence) $c(x:u)$. Consequently, so is $c^*(x:u)$ for every $x \in G \vee u$, from which the conclusion follows. ∎

## 5. Upper bound for $\omega(G \vee u)$

In propositions 5.1 through 5.5, in every instance $G \vee u$, $G$ and $\{u\}$ are necessarily vertex-disjoint. Also, $G = (V, E)$ and $G \vee u = (W, F)$ throughout this section.

**Proposition 5.1**. (Upper bound for $\omega(G)$.) *If $G \vee u$ is the secondary instance then $\omega(G \vee u) \leq c^*(G \vee u)$. (This also means the ACCN converges to a desired output.)*
**Proof.** Suppose the *ACCN* returns $c^*(G \vee u) = r$. We show no clique of $G \vee u$ can have size $r + 1$. Suppose $M = \{x_1, \ldots, x_{r+1}\}$ were a clique of $G \vee u$. Then $M \subset C(x_j)$ for each $j = 1, \ldots, r + 1$. Further, each $x_j$ has at least $r$ neighbours - call them $y_{1(j)}$ through $y_{r(j)}$ - in $A_j$ (= $C(x_j)$) such that $dy_{i(j)}(A_j) \geq r$ for $i = 1, \ldots, r$. Note that the elements of $M - \{x_j\}$ suffice to substantiate this claim. Then we have $c(x_j:u) \geq r + 1$ by the definition of $c(x_j:u)$. Consequently $c^*(x_j:u) \geq r + 1$ (for $j = 1, \ldots, r + 1$). But this leads to $c^*(G \vee u)$ ( = $r$) $\geq r + 1$, a patent contradiction. Hence $\omega(G \vee u) \leq c^*(G \vee u)$. ∎

In propositions 5.2 through 5.5, we show the *ACCN* returns $\omega(G \vee u)$ where $G$ is from any of the following graph classes: split graphs, cycles, bipartite graphs (without isolated vertex) and complete graphs.

**Proposition 5.2.** *Let $G \in \mathcal{G}_1$ where $\mathcal{G}_1$ is the class of all split graphs. Then the ACCN returns $c^*(G \vee u) = r + 1$ where $r = \omega(G)$.*
**Proof.** The following theorem is from [7].
Theorem (Hammer and Simeone). Let $d_1, \ldots, d_n$ be the degree sequence of a graph $G$ of order $n$, and assume $d_j \geq d_{j+1}$ for $j = 1, \ldots, n - 1$. Let $r$ be the largest element of the set $\{i: d_i \geq i - 1\}$. Then $G$ is a split graph if and only if $\sum_{i \text{ to } r} d_i = r(r - 1) + \sum_{r+1 \text{ to } n} d_i$.

If the degree sequence of a graph $G$ verifies the equation in the above theorem, then $\omega(G) = r$. This theorem enables one to check, in polynomial time, if $G$ is a split graph and, if so, identify $\omega(G)$.

Let $G$ be a split graph. Let $\omega(G) = r$ so that $V$ is partitioned into a maximum clique $M$ (of size $r$) and an independent set $I$ (of size, say, $p$). Write $M = \{x_1, \ldots, x_r\}$ and $I = \{y_1, \ldots, y_p\}$.

Let $x_i \in M$ be arbitrary. In $G \vee u$, we have $\{u, x_1, \ldots, x_r\} \subset C(x_i)$. If $J$ denotes the subgraph induced by $C(x_i)$ then we have $dx_j(J) \geq r$ for $j = 1, \ldots, r$ and $du(J) \geq r$. Further, if some $y \in C(x_i) \cap I$, then $dy(J) \leq r$. Consequently, $c(x:u) = r + 1$. These lead to $c^*(x:u) = r + 1$ for each $x \in M \cup \{u\}$, whence $c^*(G \vee u) = r + 1$. ∎

**Corollary 5.3.** *Given $G \in \mathcal{G}_1$, the ACCN returns $\omega(G \vee u)$.*
**Proof.** Let $G \in \mathcal{G}_1$ and $\omega(G) = r$. Then $\omega(G \vee u) = r + 1$. ∎

## 6. The worst-case time complexity of the ACCN

The worst-case time complexity of the *ACCN* is analysed using the growth-rate function Big-Oh [10, 14]. Throughout this section, by the phrase "(*) is bounded by $O(n^p)$," we will mean that there exist absolute constants $c$ and $p$ such that the computational process in the place of (*) is bounded by $cn^p$ primitive computational steps [10]. Also, by "time complexity" we will mean the worst-case one. In the following subsections, $G = (V, E)$ is the primary instance, $G \vee u = (W, F)$ is the secondary instance and $|W| = n$.

### $T_1$, the time complexity of phase 1

Adding a new vertex to the vertex set of the primary instance is done in constant time. Next, expanding the adjacency matrix of $G$ to that of $G \vee u$ bounded by $O(n^2)$ time. Hence $T_1 = O(n^2)$.

### $T_2$, the time complexity of phase 2

Obviously $|A| \leq n$. Computing $A$ is bounded by $O(n^2)$. Computing $dz(A)$ for each $z \in A$ is bounded by $O(n)$. So the execution of the inner "for" loop of phase 2 is bounded by time $O(n^2)$.

Next, the computation of $c(y:u)$ for each $y \in G \vee u$ requires $n$ readings of at most $n$ degrees in $G \vee u$, each reading followed by a logical operation of the type $dz(A) \geq k$, for each $k = 2, \ldots, n$. This computational part is bounded by $O(n^2)$. Consequently, the execution of the outer "for" loop is bounded by $O(3^2)$. Thus $T_2 = O(n^3)$.

### $T_3$, the time complexity of phase 3

The first "for" loop in phase 3 is bounded by $O(n)$ because it comprises $n$ assignments of values to $n$ variables. In the next "for" loop, for each $y \in G \vee u$, computing $c^*(y:u)$ requires $n$ readings of clique ceiling numbers of the form $c(z:u)$, each reading followed by a logical operation of the type $p \leq c(z:u)$, as $z$ runs over $G \vee u$, for each $p = 2, \ldots, n$. For each $y \in G \vee u$, these computations are bounded by $O(n^2)$. Hence $T_3 = O(n^3)$.

### $T$, the time complexity of the ACCN

From the pseudo-code presented in section 4, it is clear that phase 1, phase 2 and phase 3 run in sequence, in that order, with the algorithm terminating upon the completion of phase 3. Once the control is passed to a phase, it is not returned to any earlier phase. So the time complexity of the ACCN is $T = T_1 + T_2 + T_3 = O(n^3)$.

# 7. Comparison with three upper bound formulas

We do a comparison of the upper bound computed by the *ACCN* with those computed by three existing upper bound formulas, on the class $\mathcal{G}_1$ of split graphs.

The following three upper bound formulas (named $UB_1$ through $UB_3$) with which the *ACCN* is compared are from [2] (pages 536 and 537). In these formulas, $G \vee u$ is the secondary instance of the *ACCN*. $G \vee u$ is assumed to have $n$ vertices and $m$ edges.

($UB_1$) $\omega(G \vee u) \leq 1 + \sqrt{(2m)}$,
($UB_2$) $\omega(G \vee u) \leq 1 + \sqrt{\{2m(n-1)/n\}}$, and
($UB_3$) $\omega(G \vee u) \leq 1 + \sqrt{(2m - n + 1)}$.

For any $G \in \mathcal{G}_1$, the ACCN returns $\omega(G \vee u)$ (proved in propositions 5.2). On the other hand, the upper bounds by $UB_1$ through $UB_3$ depend on $n$ and / or $m$, and so vary. Hence these three formulas do not necessarily return $\omega(G \vee u)$.

For instance, let $G \in \mathcal{G}_1$. Suppose $G$ is split into $K_{500}$ and an independent set $I$ with $|I| = 1000$ such that $1 \leq d_y \leq 5$ for each $y \in I$ and $\sum_{y \in I} d_y = 4000$. Then $G$ has 1500 vertices and 128750 edges. So $G \vee u$ has $n = 1501$ and $m = 130250$. Also, $\omega(G \vee u) = 501$. By proposition 5.2, the *ACCN* returns $c^*(G \vee u) = 501$ (which is $\omega(G \vee u)$). On the other hand, $UB_1$, $UB_2$ and $UB_3$ return, respectively, the upper bounds 511.39, 511.22 and 509.92 for $\omega(G \vee u)$.

As another example, consider a split graph $H$ of smaller order, split into $K_5$ and an independent set $I$ with $|I| = 30$ such that $1 \leq d_y \leq 3$ for each $y \in I$ and $\sum_{y \in I} d_y = 60$. Then $H$ has 35 vertices and 70 edges. So $H \vee u$ has $n = 36$ and $m = 105$. Also, $\omega(H \vee u) = 6$. By proposition 5.2, the *ACCN* returns $c^*(H \vee u) = 6$ (which is $\omega(H \vee u)$). On the other hand, $UB_1$, $UB_2$ and $UB_3$ return, respectively, the upper bounds 15.49, 15.28 and 14.22 for $\omega(H \vee u)$.

**Merits of the *ACCN***
(i) The *ACCN* runs in $O(n^3)$ time.
(ii) No ordering of the vertices is needed (proposition 4.2). Hence the *ACCN* does not require any colouring heuristic or branching at any point.

**Limitations of the *ACCN***
(i) The *ACCN* is not designed to return any clique.
(ii) In the present version of the algorithm, the requirement of secondary instance cannot be relaxed.

# 8. An example

The primary instance is $G = (V, E)$ where $V = \{1, 2, 3, 4, 5, 6, 7\}$. The secondary instance is $G \vee u$ where $u \notin V$. The adjacency list of $G \vee u$ is:
(i) $N(1) = \{u, 2, 3, 4, 5\}$, (ii) $N(2) = \{u, 1\}$, (iii) $N(3) = \{u, 1\}$, (iv) $N(4) = \{u, 1\}$, (v) $N(5) = \{u, 1, 6, 7\}$, (vi) $N(6) = \{u, 5, 7\}$ and (vii) $N(7) = \{u, 5, 6\}$, where $N(x)$ denotes the neighbourhood (in $G \vee u$) of the vertex $x \in G \vee u$.

1) Input vertex = $u$. $C(u) = V$. For each $y \in C(u)$, $A(y,u) = C(y) \cap C(u)$. In $G \vee u$, $c(u:u)$ is computed as below:
For $y = u$, $A(y,u) = \{u, 1, 2, 3, 4, 5, 6, 7\}$ and $|N(y) \cap A(y,u)| = 7$.
For $y = 1$, $A(y,u) = \{u, 1, 2, 3, 4, 5\}$ and $|N(y) \cap A(y,u)| = 4$.
For $y = 2$, $A(y,u) = \{u, 1, 2\}$ and $|N(y) \cap A(y,u)| = 2$.

For $y = 3$, $A(y,u) = \{u, 1, 3\}$ and $|N(y) \cap A(y,u)| = 2$.
For $y = 4$, $A(y,u) = \{u, 1, 4\}$ and $|N(y) \cap A(y,u)| = 2$.
For $y = 5$, $A(y,u) = \{u, 1, 5, 6, 7\}$ and $|N(y) \cap A(y,u)| = 4$.
For $y = 6$, $A(y,u) = \{u, 5, 6, 7\}$ and $|N(y) \cap A(y,u)| = 3$.
For $y = 7$, $A(y,u) = \{u, 5, 6, 7\}$ and $|N(y) \cap A(y,u)| = 3$.
Consequently, $c(u:u) = 4$.

2) Input vertex = 1. $C(1) = \{u, 1, 2, 3, 4, 5\}$. For each $y \in C(1)$, $A(y,u) = C(y) \cap C(1)$. In $G \vee u$, $c(1:u)$ is computed as below:
For $y = u$, $A(y,u) = \{u, 1, 2, 3, 4, 5\}$ and $|N(y) \cap A(y,u)| = 5$.
For $y = 1$, $A(y,u) = \{u, 1, 2, 3, 4, 5\}$ and $|N(y) \cap A(y,u)| = 5$.
For $y = 2$, $A(y,u) = \{u, 1, 2\}$ and $|N(y) \cap A(y,u)| = 2$.
For $y = 3$, $A(y,u) = \{u, 1, 3\}$ and $|N(y) \cap A(y,u)| = 2$.
For $y = 4$, $A(y,u) = \{u, 1, 4\}$ and $|N(y) \cap A(y,u)| = 2$.
For $y = 5$, $A(y,u) = \{u, 1, 5\}$ and $|N(y) \cap A(y,u)| = 2$. Consequently, $c(1:u) = 3$.
Similarly: $c(2:u) = 3$, $c(3:u) = 3$, $c(4:u) = 3$, $c(5:u) = 4$, $c(6:u) = 4$ and $c(7:u) = 4$.

Then: $c^*(1:u) = 3$, $c^*(2:u) = 3$, $c^*(3:u) = 3$, $c^*(4:u) = 3$, $c^*(5:u) = 4$, $c^*(6:u) = 4$, $c^*(7:u) = 4$ and $c^*(u:u) = 4$. So $c^*(G \vee u) = 4$. Thus $\omega(G \vee u) \leq 4$.

## 9. Concluding remarks

We have presented a polynomial-time algorithm to compute an upper bound for $\omega(G \vee u)$ where $G = (V, E)$ is a given graph. The ACCN computes using degrees of the vertices of $G \vee u$. For each $x \in G \vee u$, the ACCN uses $d_x(G \vee u)$ or $d_x(A)$ for some appropriate $A \subset V \cup \{u\}$, and finds an upper bound for the size of any clique that contains $x$. It subsequently looks for possibilities of improving this upper bound. The ACCN does not depend on any formula but probes the closed neighbourhood of each vertex in good depth, and hence shows promising performance.

## Acknowledgement

This research was supported fully by Niels Abel Foundation, Palakkad, Kerala State, India.